\documentclass[12pt]{iopart}
\usepackage{graphicx}
\usepackage{amssymb}
\begin{document}
\title{Exact analytical approach to differential equations with variable coefficients}

\author{Mauro Bologna}

\address{ Instituto de Alta Investigaci\'{o}n, Universidad
de Tarapac\'{a}-Casilla 6-D Arica, Chile }
\ead{mauroh69@gmail.com}

\begin{abstract}
This paper shows how to build a formal analytical solution for a differential equation of arbitrary order and with variable coefficients. It proofs that the most known approximated solutions for such a problem can be derived from the analytical expression presented in the paper. The formalism can be easily extended to the infinite dimensional case such as the quantum time-dependent Hamiltonian problem.
\end{abstract}

 %\submitto{\JPA} 
\maketitle

\section{Introduction}\label{secintro}
Linear differential equations are ubiquitous in all areas of science. Very typically, a scientific problem is described by a
mathematical formalism and, very likely, eventually this implies a
differential equation, possibly linear, that is to say

\begin{equation}\label{eq1}
y^{(n)}(z)+a_{n-1}(z)y^{(n-1)}(z)+\cdots +a_0(z)y(z)=0
\end{equation}
where~$y^{(j)}$ means the~$j$th derivative of~$y(z)$. It is virtually impossible to list a complete bibliography on the topic. Its applicability ranges from quantum mechanics~\cite{fetter}, astrophysics~\cite{stev}, chemistry and stochastic processes~\cite{van,mont,west2,met} etc. just to cite text books or paper collection where the reader can find problems ending into linear differential equation with variable coefficients. Also integro-differential and non-linear differential equations can be reduced to a finite or infinite series of linear differential equations. A very famous example is the Riccati equation,~$s'(z)+s(z)^{2}=f(z)$, that, via the transformation~$s(z)=y'(z)/y(z)$, can be reduced to the second order differential equation~$y''(z)=f(z) y(z)$.

It is well known that Eq.~(\ref{eq1}) can be written as a first order equation
in the matrix formalism, i.e.

\begin{equation}\label{eq2}
\frac{d}{dz}Y(z)=Y(z)M(z)
\end{equation}
or, alternatively,

\begin{equation}\label{eq2al}
\frac{d}{dz}Y(z)=M(z)Y(z)
\end{equation}
where~$Y(z)$ is the vector~$Y(z)=\{y(z),\cdots ,y^{(n-1)}(z)\}$.
More in general the matrix~$M(z)$ describes a linear system of differential equations and may have an infinite dimension such as in the time dependent Hamiltonian problem. Indeed Eqs.~(\ref{eq2}) and (\ref{eq2al}) play a fundamental role in physics
when the matrix~$M(z)$ is the Hamiltonian of the system. In
particular Eqs.~(\ref{eq2}) and Eq.~(\ref{eq2al}) may describe either the stationary
Schr\"{o}dinger equation or the time-dependent Schr\"{o}dinger
equation, i.e. the two equations

\begin{eqnarray}\label{eq3}
&-&\frac{\hbar^2}{2m}\frac{\partial^2}{\partial x^2} \psi(x)+V(x)
\psi(x)= E \psi(x)\\\label{eq3b} & &\imath\hbar
\frac{\partial}{\partial t}\psi(x,t)=H(t) \psi(x,t)
\end{eqnarray}
The aim of this paper is to find a formal analytical solution of Eq.~(\ref{eq2}) [or Eq.~(\ref{eq2al})].
It will be shown that, starting from the analytical solution, it is possible to recover and improve the most important approximations known in literature as part of one mathematical scheme (see Ref.~\cite{kev,ben} as excellent review and text books about approximated method).

The paper is organized as follows. Section~\ref{secformsolut} shows how to build a formal
solution for a differential equation with variable coefficients. Section~\ref{lnodewvc} presents few example of exact solutions obtained using the methods of Section~\ref{secformsolut}. Section~\ref{secHamiltonian} applies the results of the
previous sections to the time dependent Hamiltonian problem. Section~\ref{secorddifequ} shows
how to recover the solution of an ordinary differential equation using the time-dependent Hamiltonian approach
previously developed. Section~\ref{sec_autovalori} shows and tests an analytical formula (derived from the previous sections) to evaluate the eigenvalues in the problem of a differential equation with boundary conditions. Finally Section~\ref{concluding} summarizes the obtained results and the future work.

\section{Formal solution for linear~$n$-order differential equations with variable coefficients}\label{secformsolut}
To solve the problem set by Eq.~(\ref{eq1}), we may solve the equivalent problem under the form of Eq.~(\ref{eq2}) [or Eq.~(\ref{eq2al})]. Our starting point is to evaluate the high order derivatives of Eq.~(\ref{eq2}). Considering the second derivative we have

\begin{equation}\label{fmatr3}
Y^{(2)}(z)=Y(z)\dot{M}(z)+ Y(z)M^2(z)=Y(z)[OM(z)]
\end{equation}
where we defined

\begin{equation}\label{o}
O\equiv \frac{d}{dz}+M(z).
\end{equation}
It is straightforward to show
that, for~$n\geq 1$

\begin{equation}\label{fmatr4}
Y^{(n+1)}(z)= Y(z)[O^n M(z)].
\end{equation}
Performing the Taylor expansion of the vectorial function~$Y(z_0+z)$ at the point~$z=z_0$, we may write

\begin{equation}\label{fmatr5}
Y(z_0+z)=
\sum_{n=0}^{\infty}Y^{(n)}(z_0)\frac{z^n}{n!}=Y(z_0)+Y(z_0)\sum_{n=0}^{\infty}\left[\frac{O^n
}{(n+1)!} M(z_0)\right]z^{n+1}.
\end{equation}
Without loss of generality we set~$z_0=0$ and, after
few algebraic steps, we can write the formal solution
for~$Y(z)$ as

\begin{equation}\label{fmatr6}
Y(z)= Y(0)\left[I+\int_0^z\exp\left[u
O\right]M(z_0)\mid_{z_0=0}du\right]
\end{equation}
where~$I$ is the identity matrix and it is understood that~$O$ acts on the right, i.e.
on~$M(z_0)$. The symbol~$\exp\left[u O\right]M(z_0)\mid_{z_0=0}$ means that it must be
evaluated at~$z_0=0$. Analogously we may write the solution for

\begin{equation}\label{fmatr6bb}
\frac{d}{dz}Y(z)=M(z)Y(z)
\end{equation}
as

\begin{equation}\label{fmatr7}
Y(z)=\left[I+\int_0^zM(z_0)\exp\left[u
O\right]\mid_{z_0=0}du\right]Y(0)
\end{equation}
where it is understood that~$O$ acts on the left.
Being the operator~$O$ the sum of two operators we can use several
tools present in literature to manipulate the solution given by
Eqs.~(\ref{fmatr6}) or (\ref{fmatr7}). Mostly there are two
important ways to proceed: the Baker-Campbell-Hausdorff
formula~\cite{baker,camp,hausd} (and its dual, the Zassenhaus
formula~\cite{magnus}) and the Trotter product formula~\cite{nelson}
and consequently the path integral technique~\cite{schulman}. The
former is left to a future work while, in this paper, we will focus on the latter.

\section{Linear~$n$-order differential equations with variable coefficients: exact examples}\label{lnodewvc}
In this section we consider examples where the formal solution gives exact results. As first example let us consider the case when~$M(z)=M_0$ i.e. a
constant matrix. Since~$[\frac{d}{dz},M_0]=0$, we may factorize
the exponential
$$\exp\left[u\left(\frac{d}{dz}+M_0\right)\right]M_0= \exp\left[u
\frac{d}{dz}\right]\exp\left[u M_0\right]M_0=\exp\left[u
M_0\right]M_0,$$ and resolving the integral of Eq.~(\ref{fmatr6}),
we obtain the well known result

\begin{equation}\label{bfmatr7}
Y(z)= Y(0) \exp\left[z M_0\right].
\end{equation}
Let us focus on the case where~$M(z)$ has polynomial elements. Polynomial coefficients in differential equations are widely studied in literature and there are many works on this topic (here a non exhaustive list of references~\cite{ref1,ref2,ref3,ref4,ref5,ref6,ref7,ref8}). We limit ourselves to few exact example in the case when the matrix elements are polynomials satisfying to certain requirements.
Setting constraints on the matrix~$M(z)$ we can find analytical
expression for Eq.~(\ref{fmatr6}). For example, requiring that

\begin{eqnarray} \label{con1}
M(z)^2=\textrm{constant}I,\,\,\,\,\frac{d^2}{dz^2}M(z)=0,
\end{eqnarray}
and using the properties of the operator~$O$, it can be showed that

\begin{eqnarray}
O^{2n}M(z)=M(z)(OM(z))^{n},\\
O^{2n+1}M(z)=(OM(z))^{n+1}.
\end{eqnarray}
Defining the matrix~$A^2(z)\equiv OM(z)$
we obtain

\begin{equation}\label{bfmatr8}
Y(z)= Y(0)\left[I+M_0\int_0^z\cosh\left[u
A_0\right]du+A_0\int_0^z\sinh\left[u A_0\right]du\right].
\end{equation}
where for compactness~$M_0\equiv M(0)$ and~$A_0\equiv A(0)$. If
$A_0$ is an invertible matrix, then expression~(\ref{bfmatr8}) may
be written

\begin{equation}\label{bfmatr9}
Y(z)= Y(0)\left[M_0A_0^{-1}\sinh\left[z A_0\right]+ \cosh\left[z
A_0\right]\right].
\end{equation}
Note that we did not specify the dimension of the matrix~$M(z)$.
This implies that solution~(\ref{bfmatr8}) holds for an arbitrary
order  of differential equations described by conditions (\ref{con1}). More: the second condition (\ref{con1}) implies that the elements of the matrix are first-degree polynomials. If we consider a~$2\times 2$ matrix, we may use~$M(z)=\sum_{i=1}^{3}m_i \sigma_i$ where~$m_i$ are assumed to be first-degree real (complex) polynomials in the real (complex) variable~$z$ and~$\sigma_i$ are the Pauli's matrices~\cite{landau3}. With this choice the condition~$M(z)^2=\textrm{constant}I$ implies the following constraint

\begin{eqnarray} \label{boh123}
M(z)^2=\sum_{i=1}^{3}m_i^2I =(a z^{2}+b z+c)I=\textrm{constant}I\Rightarrow\,\, a=0,\,\,b=0.
\end{eqnarray}
where~$a$ and~$b$ are functions of the coefficients of the polynomials~$m_i$. Since each polynomial has two coefficients we conclude that we have~$2\times 3-2=4$ degrees of freedom with respect to the choice of~$m_i$. Increasing the dimension of the matrix~$M$ it increases also the freedom's degrees. This is due to the fact that conditions~(\ref{con1}) give two constraints, independently of the dimension of the matrix~$M$. To better clarify this point, let us consider the matrix~$M(z)=\sum_{\alpha=0}^{3}m_\alpha\gamma^{\alpha}$ where~$\gamma^{\alpha}$ are the Dirac's matrices~\cite{landau4}. We have

\begin{eqnarray} \label{boh1234}
M(z)^2=\sum_{\alpha=0}^{3}m_\alpha m^\alpha I=(a z^{2}+b z+c)I=\textrm{constant}I\Rightarrow\,\, a=0,\,\,b=0.
\end{eqnarray}
As before, conditions~(\ref{con1}) give two constraints so that the arbitrary coefficients are~$2\times 4 -2=6$.

A more complex example. Let us consider the following
conditions

\begin{eqnarray} \label{conb1}
&&M(z)^2=0,\\
\label{conb2} &&
\dot{M}^2=0,\,\,\,\,\dot{M}\equiv\frac{d}{dz}M(z),
\\
\label{conb3}&& \ddot{M}=\textrm{constant},\,\,\,\,\ddot{M}\equiv\frac{d^2}{dz^2}M(z).
\end{eqnarray}
After a long but straightforward algebra we obtain

\begin{eqnarray} \label{res1}
&&O^{3n}M(z)=2^{n} M\ddot{M}^n,\\
\label{res2} && O^{3n+1}M(z)=2^{n} \dot{M}\ddot{M}^n,
\\
\label{res3} &&O^{3n+2}M(z)=2^{n}\left(\ddot{M}+M\dot{M}\right)\ddot{M}^n.
\end{eqnarray}

\begin{eqnarray} \label{res1sum}
&&f_1(z)\equiv\sum_{n=0}^{\infty}O^{3n}M\frac{u^{3n}}{3n!}=\frac{M_0}{3} e^{2^{1/3} A_0^{1/3} u}+\frac{2 M_0}{3} e^{-\frac{A_0^{1/3} u}{2^{2/3}}} \cos\left[\frac{\sqrt{3} A_0^{1/3} u}{2^{2/3}}\right],\\\nonumber
&&f_2(z)\equiv\sum_{n=0}^{\infty}O^{3n+1}M\frac{u^{3n+1}}{(3n+1)!}=\dot{M}_0\left(3\times 2^{1/3} A_0^{1/3}\right)^{-1} \\\label{res2sum}&&\left(e^{2^{1/3} A^{1/3}_0 u}-e^{-\frac{A_0^{1/3} u}{2^{2/3}}} \cos\left[\frac{\sqrt{3} A_0^{1/3} u}{2^{2/3}}\right]+\sqrt{3} e^{-\frac{A_0^{1/3} u}{2^{2/3}}} \sin\left[\frac{\sqrt{3} A_0^{1/3} u}{2^{2/3}}\right]\right),
\\
\nonumber &&f_3(z)\equiv\sum_{n=0}^{\infty}O^{3n+2}M\frac{u^{3n+2}}{(3n+2)!}=\left(\ddot{M}_0+M_0\dot{M}_0\right)\left(3\times 2^{1/3} A_0^{2/3}\right)^{-1}\\\label{res3sum}&&\left(e^{2^{1/3} A_0^{1/3} u}-e^{-\frac{A_0^{1/3} u}{2^{2/3}}} \cos\left[\frac{\sqrt{3} A_0^{1/3} u}{2^{2/3}}\right]-\sqrt{3} e^{-\frac{A_0^{1/3} u}{2^{2/3}}} \sin\left[\frac{\sqrt{3} A_0^{1/3} u}{2^{2/3}}\right]\right).
\end{eqnarray}
where~$A^3=\ddot{M}$. We finally obtain

\begin{equation}\label{bfinal}
Y(z)= Y(0)\left[I+ \int_0^z \sum_{i=1}^{3}f_i(u)du\right].
\end{equation}
Note that the integral of the~$f_i$ functions can be easily performed. As before the solution holds true for an arbitrary dimension matrix~$M$.

\section{Time dependent Hamiltonian}\label{secHamiltonian}
The solution of time dependent Hamiltonian is a fundamental topic in quantum mechanics, quantum field theory, quantum theory of many particles etc. (see for example Ref.~\cite{fetter} as excellent text and review book). Using the previous results we may write a formal solution of the quantum problem under consideration. From Eq.~(\ref{eq3b}) we infer that the operator~$O$ has the following expression

\begin{equation}\label{fmatr8_a}
O=\frac{\partial}{\partial t_0}-\imath H(t_0).
\end{equation}
Using the formal solution (\ref{fmatr7}), via Trotter formula, we may rewrite the exponential of the
operator~$O$ as

\begin{equation}\label{fmatr8}
\exp\left[u
O\right]=\lim_{N\to\infty}\prod_{j=1}^{N}\exp\left[-\imath\frac{u}{N}
 H(t_0)\right]\exp\left[\frac{u}{N}
\partial_{t_0}\right]
\end{equation}
where for simplicity we set~$\hbar=1$ and it is understood that
the operator~$\exp\left[\frac{u}{N}
\partial_{t_0}\right]$ acts on the left. Combining
Eqs.~(\ref{fmatr7}) and (\ref{fmatr8}) we can write the ket
corresponding to the solution of Eq.~(\ref{eq3b}) as

\begin{eqnarray}\nonumber
\mid\!\!\psi_t\rangle=\mid\!\!\psi_0\rangle+\\\label{fmatr9}
\lim_{N\to\infty}\prod_{j=1}^{N}\int\limits_0^t\left[-\imath
H(t_0)\right]\exp\left[-\imath\frac{u}{N}
 H(t_0)\right]\exp\left[\frac{u}{N}
\partial_{t_0}\right]\mid_{t_0=0}du\mid\!\!\psi_0\rangle.
\end{eqnarray}\nonumber
Note that the operator~$\exp\left[\frac{u}{N}
\partial_{t_0}\right]$ translates of a quantity~$\frac{u}{N}$ the
argument of the Hamiltonian~$H(t_0)$. The wave function is

\begin{eqnarray}\nonumber
\psi(x,t)=\langle x\!\! \mid\!\!\psi_t\rangle=\psi(x,0)+\\\label{fmatr10}
\lim_{N\to\infty}\int\limits_0^t\prod_{j=1}^{N}\langle x\!\!\mid
\!\!\left[-\imath
H\left(u\right)\right]\exp\left[-\imath\frac{u}{N}
H\left(j\frac{u}{N}\right)\right] du\mid\!\!\psi_0\rangle.
\end{eqnarray}
It is worthy to stress that in this formalism the time-dependence of the Hamiltonian, i.e. the dependence on~$t_0$,
is transformed in a parameter dependence. The integration variable
that gives rise to the time dependence is~$u$.

We achieve a first result. Taking the time derivative of Eq.~(\ref{fmatr9}), after little algebra, we have that for the state~$\mid\!\!\psi_t\rangle$ holds

\begin{eqnarray}\label{giral}
\mid\!\!\psi_t\rangle=
\lim_{N\to\infty} \prod_{j=1}^{N} \exp\left[-\imath\frac{t}{N}
H\left(j\frac{t}{N}\right)\right] \mid\!\!\psi_0\rangle.
\end{eqnarray}
The above expression is nothing but that the solution to the Hamiltonian problem using a discrete approach~\cite{giraldi}.

We now assume that we know
the orthonormal eigenstates of the Hamiltonian,
i.e.~$H(t_0)\mid\!\!n_{t_0}\rangle=E_n(t_0)\mid\!\!n_{t_0}\rangle$
where, for what we said, the "time"~$t_0$ has to be considered as
a parameter of a time independent Hamiltonian. Using the closure relation \cite{paul}

\begin{eqnarray}
 \sum_{n}\mid\!\!n_{t_{0}}\rangle\langle n_{t_{0}}\!\!\mid=I,\,\,
 \langle m_{t_{0}}\!\!\mid\!\!n_{t_{0}}\rangle=\delta_{nm},
\end{eqnarray}
we may rewrite Eq.~(\ref{fmatr10}) as

\begin{eqnarray}\nonumber
\psi(x,t)=\psi(x,0)+\\\nonumber
\lim_{N\to\infty}\sum_{n_1\cdots n_N
}\prod_{j=1}^{N}\int\limits_0^tdu\\\label{fmatr10ab}\langle n_{j+1}\!\!\mid\!\!n_{j}\rangle \left[-\imath
E_{n_N}(u)\right] \exp\left[-\imath\frac{u}{N}\sum_{j=0}^{N}
 E_{n_j}\left(j\frac{u}{N}\right)\right]\psi(x,0),
\end{eqnarray}
where~$\psi(x,0)=\langle x\mid\!\!\psi_0\rangle$, and
$\mid\!\!n_{j}\rangle\equiv \mid\!\!n_{j},j\frac{u}{N}\rangle$.
Inserting a new set of variable~$\mid\!\!x\rangle$ with the properties \cite{paul}
\begin{eqnarray}
 \int\limits_{-\infty}^{\infty}dx\mid\!\!x\rangle\langle x\!\!\mid=I,\,\,
 \langle x\!\!\mid\!\!y\rangle=\delta(x-y),
\end{eqnarray}
we obtain the alternative expression

\begin{eqnarray}\nonumber
\psi(x,t)=\psi(x,0)+\\\nonumber
\lim_{N\to\infty}\sum_{n_1\cdots n_N
}\prod_{j=1}^{N}\int\limits_0^tdu\int\limits_{-\infty}^{\infty} dx_j\\\label{fmatr10b}\psi^{\ast}_{n_{j+1}}\left(x_j\right)\psi_{n_{j}}\left(x{_j}\right) \left[-\imath
E_{n_N}(u)\right] \exp\left[-\imath\frac{u}{N}\sum_{j=0}^{N}
 E_{n_j}\left(j\frac{u}{N}\right)\right]\psi(x,0).
\end{eqnarray}
where~$\psi_{n_j}(x)\equiv\psi_{n_{j}}\left(x,j\frac{u}{N}\right) \equiv \langle x\!\! \mid\!\!n_{j}\rangle$. Eqs.~(\ref{fmatr10ab}) and (\ref{fmatr10b}) represent the exact solution of the starting problem, Eq.~(\ref{eq3b}).
Making the approximation
\begin{equation}
\Big |n_{j+1},(j+1)\frac{u}{N}\Big\rangle\approx \Big |n_{j+1},j\frac{u}{N}\Big\rangle+\frac{u}{N}\frac{\partial}{\partial z} \Big |n_{j+1},j\frac{u}{N}\Big\rangle,\,\,z\equiv j\frac{u}{N},
\end{equation}
then it holds

\begin{equation}
\langle n_{j+1}\!\!\mid\!\!n_{j}\rangle\approx\Big\langle n_{j+1},j\frac{u}{N}\Big |n_{j},j\frac{u}{N}\Big\rangle=\delta_{ n_{j+1}, n_{j}}.
\end{equation}
The sums on the states~$\mid\!\!n_j\rangle$ can be preformed using the delta Kronecker for two consecutive indices, and we have at first order approximation for~$N\to\infty$

\begin{eqnarray}\label{solmia}
\psi(x,t)&\approx&\psi(x,0)+\sum_{n
}\int\limits_0^tdu c_n\psi_{n}(x,u) \left[-\imath
E_{n}(u)\right] \exp\left[-\imath\int\limits_0^u
 E_{n}\left(z\right)dz\right]
\end{eqnarray}
where, by definition,~$c_n=\langle n_0\mid\!\!\psi_0\rangle$ and we used the continuous limit

\begin{eqnarray}\label{limite}
\frac{u}{N}\sum_{j=0}^{N}
 E_{n_j}\left(j\frac{u}{N}\right)\to\int\limits_0^u
 E_{n}\left(z\right)dz.
\end{eqnarray}
Via a simple integration by parts we obtain

\begin{eqnarray}\nonumber
\psi(x,t)&\approx&\sum_{n
}c_n\psi_{n}(x,u) \exp\left[-\imath\int\limits_0^u
 E_{n}\left(z\right)dz\right]+\\\label{eq41}
 &-&\sum_{n
}\int\limits_0^tdu c_n\frac{\partial}{\partial u}\psi_{n}(x,u) \exp\left[-\imath\int\limits_0^u
 E_{n}\left(z\right)dz\right].
\end{eqnarray}
The first term on the right side represents the adiabatic approximation (firstly introduced by Born and Fock~\cite{born}) while the second is the correction to the adiabatic term. Let us now evaluate the next order approximation. Keeping only the terms containing the factor~$u/N$ we may write the correction as

\begin{eqnarray}\nonumber
|\delta\psi\rangle&=&-\imath\sum_{n,k
}\int\limits_0^t du\int_{0}^{u}dz\int\limits_{-\infty}^{\infty} dy E_{n}\left(u \right)c_k\frac{\partial}{\partial z}\langle n,z\!\!\mid\!\!k,z\rangle\times
\\\label{eqcorra}
 &\times&\exp\left[-\imath\left(\int\limits_0^z
 E_{k}\left(p\right)dp+\int\limits_z^u
 E_{n}\left(p\right)dp\right)\right]\mid\!\!n,u\rangle
\end{eqnarray}
or in terms of wave function

\begin{eqnarray}\nonumber
\delta\psi(x,t)&=&-\imath\sum_{n,k
}\int\limits_0^t du\int_{0}^{u}dz\int\limits_{-\infty}^{\infty} dy E_{n}\left(u \right)c_k\psi_{k}(y,z)
\frac{\partial}{\partial z}\psi^{\dag}_{n}(y,z)\times
\\\label{eqcorr}
 &\times&\exp\left[-\imath\left(\int\limits_0^z
 E_{k}\left(p\right)dp+\int\limits_z^u
 E_{n}\left(p\right)dp\right)\right]\psi_{n}(x,u).
\end{eqnarray}
The integral on the space variable~$y$ can be expressed in terms of transition probability and we may write in more compact way

\begin{eqnarray}\nonumber
\delta\psi(x,t)&=&-\imath\sum_{n,k
}\int\limits_0^t du\int_{0}^{u}dz E_{n}\left(u \right)c_k \frac{\partial}{\partial z}P_{kn}(z)\times
\\\label{eqcorr_prob}
 &\times&\exp\left[-\imath\left(\int\limits_0^z
 E_{k}\left(p\right)dp+\int\limits_z^u
 E_{n}\left(p\right)dp\right)\right]\psi_{n}(x,u).
\end{eqnarray}
where~$P_{kn}(u)=\int\limits_{-\infty}^{\infty}\psi^{\dag}_{n}(x,u)\psi_{k}(x,u)dx$ is the transition probability from~$k$ to~$n$. Note that the diagonal terms~$k=n$ vanish identically.
In order that the above expression represents a correction to Eq.~(\ref{solmia}) it must hold

\begin{eqnarray}\label{validity}
\parallel\delta\psi(x,t)\parallel\ll \parallel \psi(x,t)\parallel
\end{eqnarray}
where~$\psi(x,t)$ is given by Eq.~(\ref{solmia}). The above inequality fixes the range of validity of our approach.
As we will see in the next section, through this formalism, we can obtain the Wentzel-Kramers-Brillouin (WKB) approximation~\cite{kev,landau3} for ordinary differential equations.

\section{Ordinary differential equations}\label{secorddifequ}
\begin{figure}[ht]
\begin{minipage}[b]{0.45\linewidth}
\centering
\includegraphics[width=.9\textwidth]{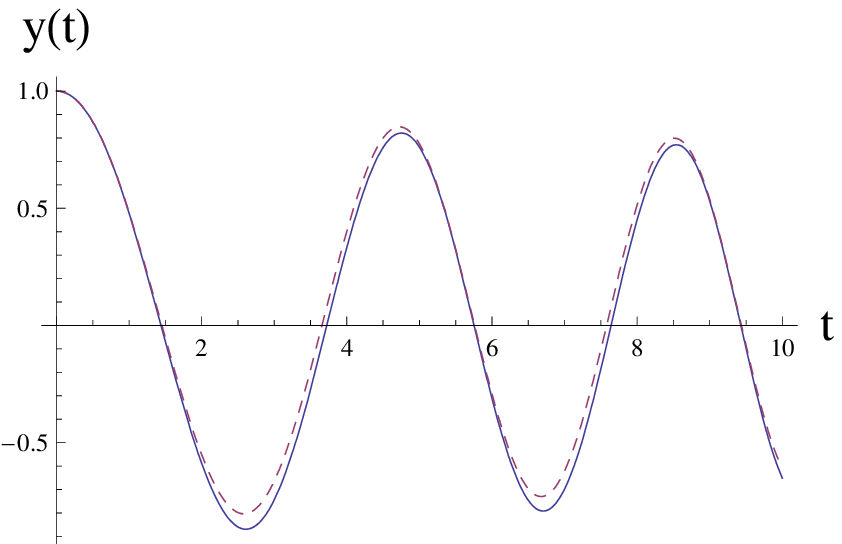}
\caption{Comparison between the exact solution (continuous line) and the approximated solution (dashed line) given by Eq.~(\ref{opla}).}
\label{nocorrected}
\end{minipage}
\hspace{0.1cm}
\begin{minipage}[b]{0.45\linewidth}
\centering
\includegraphics[width=.9\textwidth]{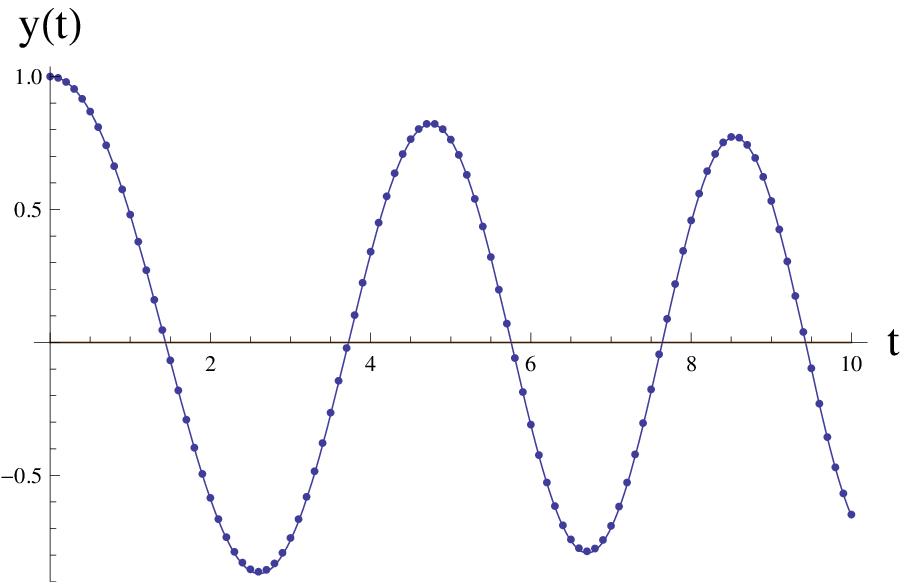}
\caption{Comparison between the exact solution (continuous line) and the approximated solution (dotted line)
given by Eq.~(\ref{opla}) where it has been added the correction given by Eq.~(\ref{eqcorra2}).}
\label{corr}
\end{minipage}
\end{figure}

With some caution, the previous approach can be extended also to ordinary differential equations,
namely when the matrix driven the system is not a hermitian operator.

\begin{equation}\label{eq2_bis}
\frac{d}{dt}Y(t)=M(t)Y(t).
\end{equation}
In the previous section an approximated solution has been obtained using the closure relation and the orthonormal condition. Since the matrix of the system is not an hermitian matrix, we need to define two sets of vectors

\begin{eqnarray}\label{normbis}
\sum_{n}\mid\!\!n\rangle\langle \bar{n}\!\!\mid=I,\,\,
 \langle \bar{m}\!\!\mid\!\!n\rangle=\delta_{nm},
\end{eqnarray}
where~$\langle \bar{n}\!\!\mid$ in general does not coincides with the hermitian conjugate of~$\mid\!\!n\rangle$. Repeating the same steps of the previous section, we end up in a similar approximated expression

\begin{eqnarray}\label{solmia_b}
\mid\!\!Y_t\rangle &\approx& \mid\!\!Y_0\rangle+\sum_{n
}\int\limits_0^tdu c_n\lambda_{n}(u)\exp\left[\int\limits_0^u
 \lambda_{n}\left(z\right)dz\right]\mid\!\!n(u)\rangle,
\end{eqnarray}
where~$c_n=\langle \bar{n}_0\!\!\mid Y_0\rangle$ and $\lambda_{n}(u)$, $\mid\!\!n(u)\rangle$ are the eigenvalues and eigenvectors of~$M(u)$ respectively. For consistency with the previous section we also introduced the symbols $\mid\!\!Y_t\rangle\equiv Y(t)$, $\mid\!\!Y_0\rangle\equiv Y(0)$. To better explain the procedure and the notation let consider the second order differential equation

\begin{equation}\label{start}
\ddot{y}=f(t)y.
\end{equation}
Writing it in matrix formalism we have

\begin{equation}\label{ex1}
\frac{d}{dt}\mid\!\!Y_t\rangle=\left[
\begin{array}{cc}
 0 & 1 \\
 f(t)& 0\\
 \end{array}
 \right]\mid\!\!Y_t\rangle,\,\,\mid Y_t\rangle=
 \left[
\begin{array}{c}
 y(t) \\
 \dot{y}(t)
 \end{array}
 \right],\,\,\mid Y_0\rangle=
 \left[
\begin{array}{c}
 a \\
 b
 \end{array}
 \right].
\end{equation}
The eigenvectors and eigenvalues of the driven matrix are

\begin{eqnarray}\label{ex1_b}
|1\rangle &=& N_1\left[
\begin{array}{c}
 1 \\
 \sqrt{f(t)}
 \end{array}
 \right],\,\,|2\rangle= N_2\left[
\begin{array}{c}
 1 \\
- \sqrt{f(t)}
 \end{array}
 \right],\\\label{ex2}\lambda_1 &=&\sqrt{f(t)},\,\, \lambda_2=-\sqrt{f(t)}
\end{eqnarray}
with~$N_i$ a factor that has to be chosen in such a way to satisfy conditions (\ref{normbis}). It is straightforward to obtain

\begin{eqnarray}\label{ex1_3}N_{1,2}&=&\frac{1}{\sqrt{2}\sqrt[4]{f(t)}}\\\label{ex3}
\langle \bar{1}|&=&\frac{1}{\sqrt{2}\sqrt[4]{f(t)}}
\left[\sqrt{f(t)}\,\,1\right],\,\,\langle \bar{2}|=\frac{1}{\sqrt{2}\sqrt[4]{f(t)}}
\left[\sqrt{f(t)}\,\,-1\right].
\end{eqnarray}
Consequently we have

\begin{eqnarray}\nonumber
y_1(t)&=&y(t)\approx a+\int_0^{t}\sqrt[4]{\frac{f(u)}{f(0)}}\times
\\\label{opla}
&&\left(a\sqrt{f(0)}\sinh\left[ \int_0^{u}\sqrt{f(z)}dz\right]+b\cosh\left[ \int_0^{u}\sqrt{f(z)}dz\right]\right)du,\\\nonumber
y_2(t)&=&\dot{y}(t)\approx b+\int_0^{t}\sqrt[4]{\frac{f(u)}{f(0)}}\sqrt{f(u)}\times\\\label{opla_bis}&&\left(a\sqrt{f(0)}\cosh\left[ \int_0^{u}\sqrt{f(z)}dz\right]+b\sinh\left[ \int_0^{u}\sqrt{f(z)}dz\right]\right)du.
\end{eqnarray}
Performing an integration by parts and redefining the constants~$\sqrt[4]{f(0)}a=A$ and~$ \left(\sqrt[4]{f(0)}\right)^{-1}b=B$ we may write Eqs.~(\ref{opla}) and (\ref{opla_bis}) in the following form

\begin{eqnarray}\nonumber
&y(t)&\approx \frac{A}{\sqrt[4]{ f(t)}}\cosh\left[ \int_0^{t}\sqrt{f(z)}dz\right]+\frac{B}{\sqrt[4]{ f(t)}}\sinh\left[ \int_0^{t}\sqrt{f(z)}dz\right]+\\\label{gopla}
&+&\!\!\!\!\frac{1}{4}\int_0^{t}\!\!\frac{\dot{f}(u)}{[f(u)]^{5/4}}\left(\!\!A\cosh\left[ \int_0^{u}\sqrt{f(z)}dz\right]\!\!+\!\!B\sinh\left[ \int_0^{u}\sqrt{f(z)}dz\right]\right)du,\\\nonumber
&\dot{y}(t)&\approx \frac{A}{\sqrt[4]{ f(t)}}\sinh\left[ \int_0^{t}\sqrt{f(z)}dz\right]+\frac{B}{\sqrt[4]{ f(t)}}\cosh\left[ \int_0^{t}\sqrt{f(z)}dz\right]+\\\label{ggopla_bis}
&-&\!\!\frac{1}{4}\int_0^{t}\!\!\frac{\dot{f}(u)}{[f(u)]^{3/4}}\left(\!\!A\sinh\left[ \int_0^{u}\sqrt{f(z)}dz\right]\!\!+\!\!B\cosh\left[ \int_0^{u}\sqrt{f(z)}dz\right]\right)du.
\end{eqnarray}
Similar expression can be found in terms of trigonometric functions, according to the sign of $f(t)$. Also, in the first term on the right side of Eq.~(\ref{gopla}) we can recognize the WKB approximation. We now consider a further example of differential equation and compare the exact solution with the approximated solution with the correction given by Eq.~(\ref{eqcorra}). Making explicit formula (\ref{eqcorra}) for the problem under consideration, we obtain

\begin{eqnarray}\nonumber
|\delta Y(t)\rangle &=&\!\!\int\limits_0^t \!\!du \lambda(u)\int_{0}^{u}\!\!dz\frac{\dot{f}(z)}{4 f(z)}
\left(\!\!\exp\left[-\int\limits_0^z
 \lambda\left(p\right)dp+\int\limits_z^u
 \lambda\left(p\right)dp\right]c_2\mid 1\rangle+\right.
\\\label{eqcorra2}
&-& \left. \exp\left[\int\limits_0^z
 \lambda\left(p\right)dp-\int\limits_z^u
 \lambda\left(p\right)dp\right]c_1\mid 2\rangle\right)
\end{eqnarray}
where~$\lambda(u)=\sqrt{f(u)}$,~$\dot{f}(z)=\partial_{z}f(z)$,~$c_n=\langle \bar{n}_{0}|Y(0)\rangle$ and the vectors~$\mid\!\! n\rangle$ are given by Eq.~(\ref{ex1_b}). Let us consider for example the following differential equation

\begin{equation}\label{xxx}
\frac{d^{2}}{dt^{2}}y+\sqrt{t+1} y=0,\,\,y(0)=1,\,\,\dot{y}(0)=0.
\end{equation}
We can compare the exact solution with the first order approximation given by Eq.~(\ref{opla}) and the approximation where we add the correction (\ref{eqcorra2}). The result shows an excellent agreement when we add the correction (see Figs. \ref{nocorrected} and \ref{corr}).

More in general it can be shown that for a function asymptotically behaving as a power law, $f(t)\sim\pm t^{\alpha}$, with respect to the differential equation (\ref{start}) the following properties hold true (for sake of brevity we omit the analytical calculations):
\begin{itemize}
\item[] i) In order that holds condition (\ref{validity}), for an asymptotic power law behaviour with $f(t)>0$, it is sufficient that $\alpha>-2$. This constraint on the power law coincides with the region of validity of WKB approximation..
\item[] ii) For $f(t)<0$ approximated solution~(\ref{opla})gives the same asymptotic results of WKB approximation plus an extra constant for $\alpha>-2$. For correction~(\ref{eqcorra2}) to balance the extra constant it is sufficient that $-2<\alpha\lesssim 1$. In this range of the parameter $\alpha$ correction~(\ref{eqcorra2}) improves the result with respect to WKB approximation.
\end{itemize}
It is worthy to stress that condition (i) is sufficient but not necessary. A simple example is given by the discontinuous function $f(t)=a\theta(1-t)$ where $\theta(z)$ is the step function. It is easy to check that formula (\ref{opla}) gives the exact solution.
We can better appreciate the novelty of our approach studying problems in finite interval containing the turning points, i.e. the points satisfying the equation~$f(\bar{t})=0$. As we will see the accuracy of the presented approximation is very good and in the next section we will show an analytical formulation for the eigenvalue problem associated to a differential equation with boundary conditions.

We will devoted the rest of the paper to this important problem. To better show this, we consider the following example where WKB approximation can not be used, at least directly. Let us study the following differential equation

\begin{figure}[ht]
\begin{minipage}[b]{0.45\linewidth}
\centering
\includegraphics[width=.9\textwidth]{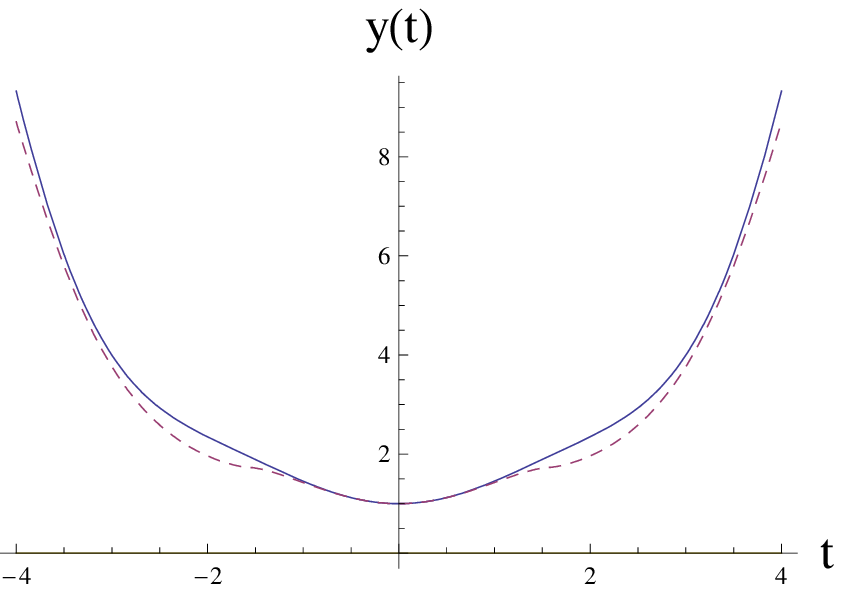}
\caption{Comparison between the exact solution (continuous line) and the approximated solution (dashed line) given by Eq.~(\ref{opla}).}
\label{figmia2}
\end{minipage}
\hspace{0.1cm}
\begin{minipage}[b]{0.45\linewidth}
\centering
\includegraphics[width=.9\textwidth]{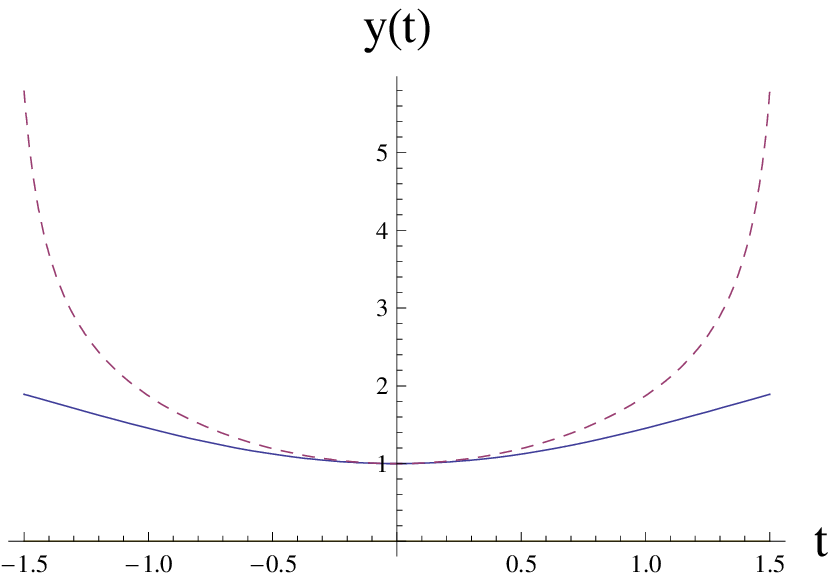}
\caption{Comparison between the exact solution (continuous line) and the approximated solution (dashed line) given by the WKB approximation. }
\label{figh2}
\end{minipage}
%\caption{Left plot: comparison between the exact result and the approximated solution given by Eq.~(\ref{solmia}). Right plot: comparison between the exact result and the approximated solution given by the adiabatic term only. Solution given by Eq.~(\ref{solmia}) more accurate}
\end{figure}
\begin{equation}\label{xx}
\frac{d^{2}}{dt^{2}}y-(\cos^{2}t) y=0,\,\,y(0)=1,\,\,\dot{y}(0)=0,\,\,\,\,t\in [-4,4].
\end{equation}
In Figs.~\ref{figmia2} and~\ref{figh2} we see that the approximation given by Eq.~(\ref{opla}) shows a small distortion with respect to the exact solution in the neighbourhood of the turning point~$t=\pi/2$ while WKB approximation, given by the formula

\begin{eqnarray}
\label{opla_no}
y(t)\approx \sqrt[4]{\frac{f(0)}{f(t)}}\cosh\left[ \int_0^{t}\sqrt{f(u)}du\right],
\end{eqnarray}
 except for a region near the origin, does not hold. As final remark of this section we note that the presented approach can be straightforwardly extended to higher order differential equations.

\section{Eigenvalue problem}\label{sec_autovalori}
In this section we will face the problem of finding the eigenvalues in the problem of a differential equation with boundary conditions. To show that, we will study in detail the following example of eigenvalue problem. Let us consider the differential equation defined in the interval $t\in [-1,1]$, that is to say

\begin{eqnarray}\frac{d^{2}}{dt^{2}}y+\lambda (1-t^{2})y=0,\,\,y(0)=0,\,\,y(1)=0,\,\,\lambda > 0.
\label{eqeig1}
\end{eqnarray}
The problem admit a general analytical solution given by

\begin{eqnarray}y(t)=e^{-\frac{1}{2} t^2 \sqrt{\lambda }}\left(a H_{\frac{1}{2} \left(-1+\sqrt{\lambda }\right)}\left(\lambda ^{1/4}t\right)+b F\left[\frac{1}{4} \left(1-\sqrt{\lambda }\right),\frac{1}{2},\sqrt{\lambda}t^2 \right]\right)
\label{eqeig2}
\end{eqnarray}
where~$H_{\nu}\left(z\right)$ is the Hermite function and~$F\left(\alpha,\beta,z\right)$ is the Kummer confluent hypergeometric function.
Imposing the condition~$y(0)=0$ we obtain

\begin{eqnarray}y(t)=A  e^{-\frac{1}{2} t^2 \sqrt{\lambda }}  t \lambda ^{1/4}  F\left[\frac{3}{4}-\frac{\sqrt{\lambda }}{4},\frac{3}{2},\sqrt{\lambda}t^2 \right],
\label{eqeig3}
\end{eqnarray}
while from the right boundary condition, i.e.~$y(1)=0$, we obtain the numerical values of~$\lambda$ that gives nontrivial solutions to Eq.~(\ref{eqeig1}). In spite of the fact that WKB approximation is diverging at~$t=1$ still it gives finite values for~$\lambda$ that vanish the function at~$t=1$

\begin{eqnarray}y_{WKB}(t)= \frac{a}{\sqrt[4]{1-t^{2}}} \sin\left[\frac{\sqrt{\lambda_{n} }}{2} \left(t \sqrt{1-t^2}+\arcsin t\right)\right],\,\,\lambda_{n}=16 n^{2}.
\label{eqeig_dirich}
\end{eqnarray}
Note that the applicability of WKB approximation rests on the fact that the divergent factor is compensated by the function sinus. Alternatively, we may evaluate the eigenvalues using the approach presented in this paper. From Eq.~(\ref{opla}), via the condition~$y(0)=0$, we have

\begin{eqnarray}
y(t)\approx a\int_0^{t}\sqrt[4]{1-u^{2}}\cos\left[\frac{\sqrt{\lambda }}{2} \left(u\sqrt{1-u^2}+\arcsin u\right)\right]du.
\label{eigen_opla}
\end{eqnarray}
Using the condition~$y(1)=0$ we have an equation for the eigenvalues. In Table \ref{table:table11} we compare the exact values of~$\lambda$ with WKB and the present approximation. The percent error is defined as $|\frac{\lambda_{ex}-\lambda_{appr}}{\lambda_{ex}}|$ where $\lambda_{ex}$ is the exact value and $\lambda_{appr}$ is the approximated value.

\begin{table}[h]
\centering
    \begin{tabular}{|l|l|l|l|l|}
    \hline
   ~$\lambda_{n}~$ (exact values) &~$\lambda_{n}~$ given by Eq.~(\ref{eigen_opla})
   &\% error &$\lambda_{n}$ given by WKB & \% error\\ \hline
  13.486& 13.767& 2\%& 16&19\% \\ \hline
    58.811& 59.174&0.6\% & 64&9\% \\ \hline
    136.140 & 136.557& 0.3\%& 144 & 6\%\\ \hline
    245.470& 245.930& 0.2\%& 256&4\%\\ \hline
   386.802 & 387.296& 0.1\%& 400& 3\% \\ \hline
    \end{tabular}\caption{A comparison between the exact values obtained from Eq.~(\ref{eqeig3}) and the approximated values given by Eq.~(\ref{opla_bis}) and WKB approximation. }
\label{table:table11}
\end{table}
In the next example we will consider a case where, due to the divergence at the turning points, WKB approximation can not be used directly. Let us modify the problem with boundary conditions stated in Eq.~(\ref{eqeig1}) as follows

\begin{eqnarray}\frac{d^{2}}{dt^{2}}y+\lambda (1-t^{2})y=0,\,\,y(x)=y(-x),\,\,y(1)=y(0),\,\,\lambda > 0.
\label{eqeig1_bis}
\end{eqnarray}
Using the approximate expression for~$y(t)$ given by Eq.~(\ref{opla}), we may write the following equation for the eigenvalues

\begin{eqnarray}
 \sqrt{\lambda}\int_{0}^{1}\left(1-u^2\right)^{1/4} \sin\left[\frac{\sqrt{\lambda }}{2} \left(u \sqrt{1-u^2 }+ \arcsin u\right)\right]du=0.
\label{eqeignewbound}
\end{eqnarray}
We obtain a sequence of of approximated eigenvalues in good agreement with the exact ones (see Table \ref{table:table3bis}). Note that the boundary condition is at the turning point, i.e.~$t=1$, and we can not use directly the WKB approximation since the approximated function is divergent at the turning point.

\begin{table}[h]
\centering
    \begin{tabular}{|l|l|l|}
    \hline
   ~$\lambda_{n}~$ (exact values) &~$\lambda_{n}~$ given by Eq.~(\ref{eqeignewbound})
   & Relative error\\ \hline
   43.185  & 46.138 & \,7\% \\ \hline
   77.736 &74.721&\,\,4\% \\ \hline
    208.573  &  213.915&\,\,2.6\%\\ \hline
    286.144 &281.010&\,\,2\% \\\hline
    500.880 &  508.297&\,\,1.5\%\\ \hline
    \end{tabular}\caption{A comparison between the exact values $\lambda_n$ and the approximated values given by Eq.~(\ref{eqeignewbound})}
\label{table:table3bis}
\end{table}
Finally we study a differential equation with mixed boundary conditions, i.e. Dirichlet-von Neumann conditions

 \begin{eqnarray}\frac{d^{2}}{dt^{2}}y+\lambda (1-t^{2})y=0,\,\,y(0)=0,\,\,\frac{d}{dt}y_{\mid_1}\equiv\dot{y}(1)=0,\,\,\lambda > 0.
\label{eqeig1_bis}
\end{eqnarray}
 Evaluating the solution at~$t=0$ we obtain again expression (\ref{eqeig3}). Vanishing the derivative of~$y(t)$ at~$t=1$ we obtain the eigenvalues of the problem. Since the condition on the derivative is at the turning point, i.e.~$t=1$, we can not use directly the WKB approximation. On the other hand we may use the approximate expression for the derivative of~$y(t)$ given by Eq.~(\ref{opla_bis}) that yields

\begin{eqnarray}\label{opla_ter}
\dot{y}(t)\approx b\left(1-\sqrt{\lambda}\int_0^{t} (1-u^{2})^{\frac{3}{4}}
\sin\left[\frac{\sqrt{\lambda }}{2} \left(u\sqrt{1-u^2}+\arcsin u\right)\right]du\right)\!\!.
\end{eqnarray}
The values of~$\lambda$ so obtained are compared with the exact ones in Table \ref{table:table2}. In general we can use Eqs.~(\ref{opla}) and (\ref{opla_bis}) as formulas to find, analytically, the eigenvalues of a differential equation with boundary conditions defined in a finite interval.

\begin{table}[h]
\centering
    \begin{tabular}{|l|l|l|}
    \hline
   ~$\lambda_{n}~$ (exact values) &~$\lambda_{n}~$ given by Eq.~(\ref{opla_ter})
   & Relative error\\ \hline
   5.122  & 4.721 & \,\,8\% \\ \hline
    39.661 & 39.836& 0.4\% \\ \hline
    106.249 & 106.063& 0.2\%\\ \hline
    204.856 &204.952& 0.05\% \\\hline
    335.473 & 335.352& 0.03\%\\ \hline
    \end{tabular}\caption{A comparison between the exact values obtained vanishing the derivative of Eq.~(\ref{eqeig3}) at $t=1$ and the approximated values given by Eq.~(\ref{opla_ter})}
\label{table:table2}
\end{table}

\section{Concluding Remarks}\label{concluding}
In this paper it has been showed a formal analytical solution for a linear differential equation of arbitrary order~$n$ and with variable coefficients. The solution has been found exploiting the reduction of the problem to a first order differential equation through the matrix formalism. This approach provided a solution also for differential equation involving operators such as the time-dependent Schr\"{o}dinger
equation.

It has been shown that the most important approximations available in literature can be deduced from the proposed solution. With respect to ordinary differential equation, a noteworthy fact is that the proposed approximation is smooth around the turning points and this allows to study solutions of differential equation at points where the WKB approximation, in general, can not be used. Examples of analytical eigenvalue evaluation have been showed. Even though we studied examples of second order differential equations, the approach can be easily extended to higher order differential equations. Finally this work showed that there are several possible research lines that have to be fully explored. In particular the approach based on Baker-Campbell-Hausdorff formula is left for a future work.

\section*{Acknowledgments}
The author acknowledges financial support from UTA Mayor project no xxx

\end{document}